\documentclass[12pt,a4paper]{article}
\usepackage{amsmath,amsfonts,amssymb,amsthm,amscd}
\newtheorem{thm}{Theorem}
\newtheorem{lem}[thm]{Lemma}
\newtheorem{prop}[thm]{Proposition}
\newtheorem{defn}[thm]{Definition}
\newtheorem{cor}[thm]{Corollary}
\newtheorem{notation}[thm]{Notations}
\newtheorem{rem}[thm]{Remark}
\begin{document}

\title{About the geometry of almost para-quaternionic
manifolds}

\author{LIANA DAVID}

\maketitle

\textbf{Abstract:} We provide a general criteria for the
integrability  of the almost para-quaternionic structure of an
almost para-quaternionic manifold $(M, {\mathcal P})$ of dimension
$4m\geq 8$ in terms of the integrability of two or three sections
of the defining rank three vector bundle $\mathcal P$. We relate
it with the integrability of the canonical almost complex
structure of the twistor space and with the integrability of the
canonical almost para-complex structure of  the reflector space of
$(M, {\mathcal P})$. We deduce that $(M, {\mathcal P})$ has plenty
of locally defined, compatible, complex
and para-complex structures, provided that $\mathcal P$ is integrable.\\

{\bf Key words and phrases:} almost para-quaternionic manifolds,
compatible complex and para-complex structures, twistor and
reflector
spaces.\\

{\bf Mathematics Subject Classification:} 53C15, 53C25.
\footnote{{\it Affiliation:} Institute of Mathematics of the
Romanian Academy "Simion Stoilow", Calea Grivitei nr. 21,
Bucharest, Romania; e-mail: liana.david@imar.ro}

\section{Introduction and main results}

An almost para-complex structure on a smooth manifold $M$  is an
endomorphism $J$ of $TM$ such that $J^{2}=  \mathrm{Id}$ (where
"$\mathrm{Id}$" denotes the identity endomorphism) and the two
distributions $\mathrm{Ker}(J\pm \mathrm{Id})$ have the same rank.
We say that $J$ is integrable (or is a para-complex structure) if
its Nijenhuis tensor
$$
N_{J}(X, Y):= [X, Y]+[JX, JY]- J\left( [JX, Y]+[X,
JY]\right),\quad\forall X, Y\in {\mathcal X}(M)
$$
is zero.

An almost para-quaternionic structure on a smooth manifold $M$ of
dimension $4m\geq 8$ is a rank three sub-bundle $\mathcal P\subset
\mathrm{End}(TM)$ locally spanned by almost para-hypercomplex
structures, i.e. by triples $\{ J_{1}, J_{2}, J_{3}\}$ where
$J_{1}$ is an almost complex structure, $J_{2}$ and $J_{3}$ are
anti-commuting almost para-complex structures and $J_{3} =
J_{1}J_{2}.$ We shall often refer to such a triple as an
admissible basis of $\mathcal P$. The bundle $\mathcal P$ comes
with a standard Lorenzian metric $\langle \cdot ,\cdot \rangle$
defined in terms of an admissible basis by
$$
\langle \sum_{i=1}^{3} a_{i}J_{i}, \sum_{j=1}^{3}b_{j}J_{j}\rangle
:=-a_{1}b_{1}+a_{2}b_{2}+a_{3}b_{3}.
$$
A para-quaternionic connection on $(M, {\mathcal P})$ is a linear
connection on $M$ which preserves the bundle $\mathcal P$. We say
that $\mathcal P$ is a para-quaternionic structure if $(M,
{\mathcal P})$ has a torsion-free para-quaternionic connection;
equivalently, if the torsion tensor of $\mathcal P$ is zero.

We begin the paper by recalling briefly, in Section \ref{Gstr},
the theory of $G$-structures. We then apply these considerations
to almost para-quaternionic manifolds. A central role in our paper
is played by the torsion tensor of an almost para-quaternionic
manifold $(M, {\mathcal P})$. In Section \ref{obatas} we give an
account on the torsion tensor of $(M, {\mathcal P})$, providing
more insight of some of the results developed in Section 2 of
\cite{zam}. We define a canonical family of para-quaternionic
connections (also called "minimal para-quaternionic connections")
on $(M,{\mathcal P})$, which consists of all para-quaternionic
connections whose torsion is equal to the torsion tensor of
$\mathcal P$. These connections are parametrized by $1$-forms and
are similar to the so called "Oproiu connections" of almost
quaternionic manifolds, which were introduced  for the first time
in \cite{oproiu} and have been used in \cite{ponte} to define a
canonical almost complex structure on the twistor space of an
almost quaternionic manifold.

In Section \ref{nou} we prove that if an almost para-quaternionic
manifold $(M, {\mathcal P})$ admits two independent, compatible,
globally defined, (para)-complex structures $I_{1}$ and $I_{2}$,
such that either $I_{1}$ or $I_{2}$ is a complex structure, or
otherwise both $I_{1}$ and $I_{2}$ are para-complex structures
and, for any $p\in M$,  $\mathrm{Span}\{I_{1}(p), I_{2}(p)\}$ is a
non-degenerate $2$-plane of ${\mathcal P}_{p}$ (with its Lorenzian
metric $\langle \cdot ,\cdot \rangle$), then $(M, {\mathcal P})$
is a para-quaternionic manifold (see Theorem \ref{sug}). (An
almost (para)-complex structure on $(M, {\mathcal P})$ is
compatible with $\mathcal P$ if it is a section of the bundle
$\mathcal P$; two almost (para)-complex structures $I_{i}$ and
$I_{j}$ are independent if $I_{i}(p)\neq \pm I_{j}(p)$ at any
point $p$). If, however, $I_{1}$ and $I_{2}$ are both para-complex
and, for any $p\in M$, the $2$-plane $\mathrm{Span}\{ I_{1}(p),
I_{2}(p)\}$ is degenerate, then we need an additional compatible
para-complex structure $I_{3}$, such that $\{ I_{1}, I_{2},
I_{3}\}$ are pairwise independent but dependent at any point, to
conclude that $\mathcal P$ is para-quaternionic (see Theorem
\ref{sug}). At the end of Section \ref{nou} we construct a class
of almost para-quaternionic manifolds $(M, {\mathcal P})$ which
are not para-quaternionic but admit three independent, globally
defined, compatible para-complex structures $I_{1}$, $I_{2}$ and
$I_{3}$, such that, at any point $p\in M$ and for any $i\neq j$,
the $2$-plane $\mathrm{Span}\{I_{i}(p), I_{j}(p)\}$ is degenerate
(see Proposition \ref{propo}). Recall that for almost quaternionic
manifolds the existence of two independent, globally defined,
compatible, complex structures insures the integrability of the
almost quaternionic structure (see Theorem 2.4 of \cite{ponte}).
For conformal oriented $4$-manifolds, the existence of three
pairwise independent, globally defined, orthogonal complex
structures is needed to deduce that the conformal structure is
self-dual (see \cite{sal}, page 121).

In Section \ref{twistor} we consider the twistor space $Z^{-}$ and
the reflector space $Z^{+}$ of $(M ,{\mathcal P})$, consisting of
all compatible, complex, respectively para-complex structures of
tangent spaces of $M$, i.e.
\begin{equation}\label{zpm}
Z^{\pm} = \{ A\in {\mathcal P}:\quad A^{2}= \pm \mathrm{Id}\}
\subset\mathrm{End}(TM).
\end{equation}
It is known that a para-quaternionic connection $\nabla$ on $(M,
{\mathcal P})$ induces an almost complex structure ${\mathcal
J}^{\nabla , -}$ (respectively, an almost para-complex structure
${\mathcal J}^{\nabla , +}$) on the twistor space $Z^{-}$
(respectively, on the reflector space $Z^{+}$) and the way
${\mathcal J}^{\nabla ,\pm}$  depend on $\nabla$ has been studied
in \cite{ivanov}. Our main observation in this setting is that
${\mathcal J}^{\nabla , \pm}$ are independent of the choice of
$\nabla$, provided that $\nabla$ is minimal. Using  minimal
para-quaternionic connections we define an almost complex
structure $\mathcal J^{-}$ on $Z^{-}$ and an almost para-complex
structure ${\mathcal J}^{+}$ on $Z^{+}$, both ${\mathcal J}^{-}$
and ${\mathcal J}^{+}$ being canonical (since they depend only on
the torsion tensor of $\mathcal P$). We use ${\mathcal J}^{\pm}$
to prove that $(M, {\mathcal P})$ is integrable if and only if it
has plenty of locally defined, compatible, complex and
para-complex structures (see Theorem \ref{final}). Similar
considerations hold for almost quaternionic manifolds, the role of
minimal connections on $(M, {\mathcal P})$ being played by the
Oproiu connections of an almost quaternionic manifold (see
\cite{ponte}). The geometry of twistor and reflector spaces of
para-quaternionic manifolds with an additional compatible metric
(the so called "para-quaternionic Hermitian" and
"para-quaternionic K\"{a}hler manifolds") has already been studied
in the literature, see for example, \cite{dav}, \cite{garcia}, \cite{ivanov1}.\\

{\bf Acknowledgements:} I am grateful to Dmitri Alekseevsky for
his encouragements and interest in this work and to Paul Gauduchon
for making useful and pertinent comments on a preliminary version
of this paper.

\section{$G$-structures}\label{Gstr}

In this Section we recall the definition of the torsion tensor of
a $G$-structure \cite{2}, \cite{4}. We follow closely the
treatment developed in \cite{gaud}, Section 2.1.

Let $G$ be a closed subgroup of the general linear group
$GL_{n}(V)$, where $V=\mathbb{R}^{n}$. A $G$-structure on an
$n$-dimensional manifold $M$ is a principal $G$ sub-bundle $P$ of
the frame bundle of $M$. A linear connection on $M$ is adapted to
the $G$-structure if it is induced by a $G$-invariant connection
on $P$. Any two adapted connections $\nabla$ and $\nabla^{\prime}$
are related by $\nabla^{\prime} = \nabla +\eta$, where $\eta\in
\Omega^{1}(M, \mathrm{ad}(P))$ is a $1$-form with values in
$\mathrm{ad}(P)$, the vector bundle on $M$ associated to the
adjoint representation of $G$ on its Lie algebra. Define the
linear torsion map
\begin{equation}\label{delta}
\delta :\Omega^{1}(M,\mathrm{ad}(P))\to \Omega^{2}(M, TM),\quad (
\delta \eta )(X, Y):=\eta (X,Y)-\eta (Y,X),
\end{equation}
where $X, Y\in TM$. The image of the torsion
$T^{\nabla}\in\Omega^{2}(M, TM)$ of an adapted connection $\nabla$
into the quotient space $\frac{\Omega^{2}(M,
TM)}{\mathrm{Im}\delta}$ is independent of the choice of $\nabla$
and is called the torsion tensor of the $G$-structure $P$. It will
be denoted by $T^{P}.$

Suppose now that it is given a complement $C(\mathrm{ad}(P))$ of
$\delta \Omega^{1}(M, \mathrm{ad}(P))$ in $\Omega^{2}(M, TM)$,
i.e. a decomposition
\begin{equation}\label{dec}
\Omega^{2}(M, TM) = \delta\Omega^{1}(M, \mathrm{ad}(P)) \oplus
C(\mathrm{ad}(P)).
\end{equation}
The decomposition (\ref{dec}) identifies the quotient
$\frac{\Omega^{2}(M, TM)}{\mathrm{Im}\delta}$ with
$C(\mathrm{ad}(P)).$ An adapted connection with torsion equal to
$T^{P}\in C(\mathrm{ad}(P))$ is called minimal. Any two minimal
connections $\nabla$ and $\nabla^{\prime}$ are related by
$\nabla^{\prime} =\nabla + \eta$, where $\eta\in \Omega^{1}(M,
\mathrm{ad}(P))$ belongs to the kernel of $\delta .$

\section{Almost para-quaternionic manifolds}\label{obatas}

\subsection{Torsion of almost para-quaternionic manifolds}

Let $(M, {\mathcal P})$ be an almost para-quaternionic manifold of
dimension $n=4m\geq 8$ (in this paper we will always assume that
the dimension of the almost para-quaternionic manifolds is bigger
or equal to eight). The almost para-quaternionic structure
$\mathcal P$ defines a $G = GL_{m}(\mathbb{H}^{+})Sp(1,
\mathbb{R})$ structure on $M$, where $Sp(1, \mathbb{R})$ is the
group of unit para-quaternions acting on $\mathbb{R}^{n}$ and
$GL_{m}(\mathbb{H}^{+})$ is the group of automorphisms which
commutes with the action of $Sp(1,\mathbb{R})$ (for details see,
for example, \cite{ivanov}). We denote by $Z({\mathcal P})$ and
$N({\mathcal P})= Z({\mathcal P})\oplus {\mathcal P}$ the
centralizer, respectively the normalizer of ${\mathcal P} $ in
$\mathrm{End}(TM)$. They are vector bundles on $M$ associated to
the adjoint representations of $GL_{m}({\mathbb{H}}^{+})$ and $G$
on their Lie algebras.

The aim of this Section is to show that $\delta\Omega^{1}(M,
Z({\mathcal P}))$ and $\delta \Omega^{1}(M, N({\mathcal P}))$ have
canonical complements in $\Omega^{2}(M, TM)$, where $\delta$ is
the linear torsion map. We then relate the torsion tensor
$T^{\mathcal P}$ of $\mathcal P$ with the torsion tensor $T^{H}$
of any compatible almost para-hypercomplex structure $H$ and we
determine conditions on $T^{H}$ which insure the integrability of
$\mathcal P .$ We shall need these considerations (especially
Corollary \ref{fol}) in the proof of Theorem \ref{sug}. Our
arguments are similar to those employed in \cite{gaud} and
\cite{ponte}. This Section is intended mostly for completeness of
the text: except the different treatment, some results of this
Section were already proved in \cite{zam}.

\begin{notation}\label{epsilonc}{\rm To unify notations, we define an
(almost) $\epsilon$-complex structure on $M$ (with $\epsilon = \pm
1$) to be an (almost) complex structure when $\epsilon =-1$ and an
(almost) para-complex structure when $\epsilon =1$. In particular,
the Nijenuis tensor of an almost $\epsilon$-complex structure $J$
is
$$
N_{J}(X, Y) = \epsilon [X, Y] +[JX, JY] - J\left( [JX, Y] + [X,
JY]\right) ,\quad\forall X, Y\in {\mathcal X}(M).
$$
For an admissible basis $\{ J_{1}, J_{2}, J_{3}\}$ of $\mathcal
P$, we define $\epsilon_{i}\in \{ -1, +1\}$ by the conditions
$J_{i}^{2} = \epsilon_{i}\mathrm{Id}$, for any $i\in \{ 1,2,3\}$;
hence $\epsilon_{1}=-1$ and $\epsilon_{2}=\epsilon_{3}=1$.}
\end{notation}

\begin{notation}{\rm An operator, expression, etc,
defined in terms of an admissible basis of $\mathcal P$ but
independent of the choice of admissible basis will be considered,
without further explanation, defined on $M$.}
\end{notation}

In the next Lemma we show  that $\delta (\Omega^{1}(M, Z({\mathcal
P}))$ has a canonical complement in $\Omega^{2}(M, TM)$.

\begin{lem}\label{h1} Let $\{ J_{1}, J_{2}, J_{3}\} $ be an admissible basis of
$\mathcal P .$ Define an endomorphism
\begin{equation}\label{p}
P:\Omega^{2}(M, TM)\to \Omega^{2}(M, TM),\quad
P(T):=\frac{2}{3}\sum_{i=1}^{3} \Pi^{0,2}_{J_{i}}(T) ,
\end{equation}
where, for any $T \in \Omega^{2}(M, TM)$ and $X, Y\in TM$,
$$
\Pi^{0,2}_{J_{i}}(T )(X, Y):=\frac{1}{4}\{ T (X, Y)+\epsilon_{i}T
(J_{i}X, J_{i}Y)-\epsilon_{i}J_{i} \left( T (J_{i}X, Y)+T (X,
J_{i}Y)\right)\} .
$$
Then $P$ is independent of the choice of $\{ J_{1}, J_{2},
J_{3}\}$, is a projector (i.e. $P^{2} = P$) and $\mathrm{Ker}(P )=
\delta \Omega^{1}(M, Z({\mathcal P})).$ In particular,
\begin{equation}\label{d}
\Omega^{2}(M, TM) = \delta \Omega^{1}(M, Z({\mathcal P})) \oplus
\mathrm{Im}(P).
\end{equation}

\end{lem}

\begin{proof} Note that the expressions
$\sum_{i=1}^{3}\epsilon_{i}T (J_{i}X, J_{i}Y)$,
$\sum_{i=1}^{3}\epsilon_{i}J_{i}T (J_{i}X, Y)$ and
$\sum_{i=1}^{3}\epsilon_{i}J_{i}T (X, J_{i}Y)$ are independent of
the choice of admissible basis of $\mathcal P$. The same holds for
$P$, which is a well-defined endomorphism of $\Omega^{2}(M, TM).$

We now prove that $P$ is a projector. We shall use the notation
$J_{ij}$ for the composition $J_{i}\circ J_{j}$ and
$\epsilon_{ij}$ for $\epsilon_{i}\epsilon_{j}$. For any $i$,
$\Pi^{0,2}_{J_{i}}$ is a projector of $\Omega^{2}(M, TM)$ and for
any $i\neq j$,
$$
\Pi^{0,2}_{J_{i}}\circ \Pi^{0,2}_{J_{j}} = \frac{1}{4} \left(
\Pi^{0,2}_{J_{i}}+\Pi^{0,2}_{J_{j}} -\Pi^{0,2}_{J_{ij}}\right)
+\frac{1}{16}{\mathcal E}_{ij},
$$
where the endomorphism ${\mathcal E}_{ij}$ of $\Omega^{2}(M, TM)$
has the following expression: for any $T\in \Omega^{2}(M, TM)$ and
$X, Y\in TM$,
\begin{align*}
{\mathcal E}_{ij}(T )(X, Y)& = \epsilon_{ij}J_{j}\left( T
(J_{ij}X, J_{i}Y) +T (J_{i}X, J_{ij}Y)\right)\\
&+\epsilon_{ij}J_{i}\left( T (J_{ij}X, J_{j}Y) +T
(J_{j}X, J_{ij}Y)\right)\\
&+\epsilon_{ij}J_{ij}\left( T (J_{i}X, J_{j}Y) +T (J_{j}X,
J_{i}Y)\right).
\end{align*}
Since ${\mathcal E}_{ij}$ is anti-symmetric in $i$ and $j$,
\begin{equation}\label{pi}
\Pi^{0,2}_{J_{i}}\circ \Pi^{0,2}_{J_{j}}+ \Pi^{0,2}_{J_{j}}\circ
\Pi^{0,2}_{J_{i}} = \frac{1}{2}\left(
\Pi^{0,2}_{J_{i}}+\Pi^{0,2}_{J_{j}}- \Pi^{0,2}_{J_{ij}}\right) ,
\end{equation}
for any $i\neq j.$ Relation (\ref{pi}) implies that $P^{2}=P$. We
now prove that $\mathrm{Ker} (P)= \delta\Omega^{1}(M, Z({\mathcal
P})).$ It is easy to show, using definitions,  that $\delta
\Omega^{1}(M,Z({\mathcal P}))$ is included in the kernel of $P$.
Conversely, in order to show that $\mathrm{Ker}(P)\subset \delta
\Omega^{1}(M,Z({\mathcal P}))$ we define an endomorphism $\pi$ of
$\Omega^{2}(M, TM)$ by
\begin{align*}
\pi (T)(X, Y)& := \frac{1}{4}T(X, Y)+\frac{1}{4}\sum_{i=1}^{3}\epsilon_{i}J_{i}T(X, J_{i}Y)\\
&-\frac{1}{12}\sum_{i=1}^{3}\epsilon_{i}  J_{i}T(J_{i}X,
Y)-\frac{1}{12}\sum_{i, j=1}^{3}\epsilon_{i}\epsilon_{j}
J_{j}J_{i}T(J_{i}X, J_{j}Y),
\end{align*}
where $T\in \Omega^{2}(M,TM)$ and $X, Y\in TM.$ The endomorphism
$\pi $ is independent of the choice of admissible basis of
$\mathcal P$ and its image is included in $\Omega^{1}(M,
Z({\mathcal P}))$. Moreover, it can be checked that
\begin{equation}\label{pi1}
\delta\circ \pi (T)=T-P(T),\quad\forall T\in \Omega^{2}(M, TM).
\end{equation}
Relation (\ref{pi1}) implies the converse inclusion
$\mathrm{Ker}(P)\subset \delta \Omega^{1}(M, Z({\mathcal P}))$. We
deduce that $\mathrm{Ker}(P)= \delta \Omega^{1}(M, Z({\mathcal
P}))$. Since $P^{2}=P$, the decomposition (\ref{d}) follows.

\end{proof}

\begin{cor}\label{torsionc}
Let $H= \{ J_{1}, J_{2}, J_{3}\}$ be an admissible basis of
$\mathcal P$ and $\nabla$ a linear connection which preserves all
$J_{i}$. Then
\begin{equation}\label{ptn}
T^{H}:=P(T^{\nabla})= -\frac{1}{6}\sum_{i=1}^{3}
\epsilon_{i}N_{J_{i}}
\end{equation}
is independent of the choice of $\nabla$ and is the torsion tensor
of the almost para-hypercomplex structure $H$.
\end{cor}

\begin{proof}
If $J$ is an almost $\epsilon$-complex structure on a manifold $M$
and $\nabla$ is a linear connection on $M$ which preserves $J$,
then
\begin{equation}\label{N}
\Pi^{0,2}_{J}(T^{\nabla})(X, Y) = -\frac{\epsilon}{4} N_{J}(X, Y),
\quad\forall X, Y\in TM.
\end{equation}
The first claim follows from (\ref{N}) and the definition of $P$.
The second claim is trivial, from  our considerations of Section
\ref{Gstr} and from Lemma \ref{h1}.

\end{proof}

\begin{rem}{\rm The linear torsion map $\delta$ is injective for
$GL_{m}(\mathbb{H}^{+})$-structures. Given an almost
para-hypercomplex structure $H$ there is a unique linear
connection, called the Obata connection, which preserves $H$ and
whose torsion is equal to $T^{H}$ (see \cite{zam}, Proposition
2.1).}
\end{rem}

We need the following Lemma for the proof of Proposition
\ref{decn}.

\begin{lem}\label{ajut}
Let $P$ be the projector defined in Lemma \ref{h1}. For any
admissible basis $\{ J_{1}, J_{2}, J_{3}\}$ of $\mathcal P$ and
$T\in \mathrm{Im}(P)\subset\Omega^{2}(M, TM)$,
\begin{equation}\label{d1}
\sum_{i=1}^{3}\epsilon_{i} J_{i} T(J_{i}X, Y) = - T(X,
Y),\quad\forall X, Y\in TM.
\end{equation}
In particular,
\begin{equation}\label{d2}
\sum_{i=1}^{3} \epsilon_{i} \mathrm{tr} \left( J_{i}
T_{J_{i}X}\right)  = 0,\quad\forall X \in TM.
\end{equation}
Above, $T_{Y}:= T(Y, \cdot )$ denotes the interior product of a
tangent vector $Y\in TM$ with the $TM$-valued $2$-form $T.$

\end{lem}

\begin{proof} Relation (\ref{d1}) can be checked by
writing $T = P(A)$, for $A\in \Omega^{2}(M, TM)$, and using the
definition of $P$. Relation (\ref{d2}) follows from (\ref{d1}) and
the observation that $\mathrm{tr}\left(T_{X} \right)=0$ for any
$T\in \mathrm{Im}(P)$ and $X\in TM.$
\end{proof}

We now state the main result of this Section.

\begin{prop}\label{decn}
Let $\{ J_{1}, J_{2}, J_{3}\}$ be an admissible basis of $\mathcal
P .$ The subspace
$$
C(N({\mathcal P})):= \{ T\in \mathrm{Im}(P):\quad
\mathrm{tr}\left( J_{i}T_{X}\right) =0,\quad\forall X\in
TM,\quad\forall i=1,2,3\}
$$
is a complement of $\delta \Omega^{1}( M, N({\mathcal P}))$ in
$\Omega^{2}(M, TM).$ The projection on the second factor of the
decomposition
\begin{equation}\label{dnq}
\Omega^{2}(M, TM)= \delta \Omega^{1}( M, N({\mathcal P}))\oplus
C(N({\mathcal P}))
\end{equation}
is the map
\begin{equation}\label{proj}
\Omega^{2}(M,TM)\ni T\rightarrow P(T)-\delta (\tau_{1}\otimes
J_{1}+\tau_{2}\otimes J_{2}+\tau_{3}\otimes J_{3}),
\end{equation}
where $P$ is the projector of Lemma \ref{h1} and the $1$-forms
$\tau_{i}$ are defined by
\begin{equation}\label{taui}
\tau_{i}(X)=\frac{\epsilon_{i}\mathrm{tr}\left(
J_{i}P(T)_{X}\right)}{n-2},\quad  \forall i\in \{ 1,2,3\} .
\end{equation}
\end{prop}

\begin{proof} We define two subspaces of $\Omega^{1}(M, {\mathcal
P})$:
\begin{align*}
\Omega^{1}_{0}(M, {\mathcal P})&:=
\{\sum_{i=1}^{3}\gamma_{i}\otimes J_{i}:\quad
\sum_{i=1}^{3}\gamma_{i}(J_{i}X)=0,\quad\forall  X\in TM\};\\
{\mathcal C}(M)&:= \{ T^{\alpha}:=\sum_{i=1}^{3}\epsilon_{i}
(\alpha\circ J_{i})\otimes J_{i},\quad \forall\alpha\in
\Omega^{1}(M)\} .
\end{align*}
Clearly, $\Omega^{1}_{0}(M, {\mathcal P})$ and ${\mathcal C}(M)$
have trivial intersection. Consider now an arbitrary $\mathcal
P$-valued $1$-form $\sum_{i=1}^{3}\alpha_{i}\otimes J_{i}$ and
define $\alpha := -\frac{1}{3}\sum_{i=1}^{3}\alpha_{i}\circ
J_{i}$. Then
\begin{align*}
A:= \sum_{i=1}^{3}\alpha_{i}\otimes J_{i} + T^{\alpha}
&=\left(\frac{2}{3}\alpha_{1}-\frac{1}{3}\alpha_{2}\circ
J_{3}+\frac{1}{3}\alpha_{3}\circ J_{2}\right)\otimes J_{1}\\
&+\left(\frac{2}{3}\alpha_{2}-\frac{1}{3}\alpha_{1}\circ
J_{3}-\frac{1}{3}\alpha_{3}\circ J_{1}\right)\otimes J_{2}\\
&+\left(\frac{2}{3}\alpha_{3} +\frac{1}{3}\alpha_{1}\circ
J_{2}+\frac{1}{3}\alpha_{2}\circ J_{1}\right)\otimes J_{3}
\end{align*}
belongs to $\Omega^{2}_{0}(M, {\mathcal P}).$ Therefore,
$\Omega^{1}(M, {\mathcal P})$ decomposes as
\begin{equation}\label{initial}
\Omega^{1}(M, {\mathcal P}) ={\mathcal
C}(M)\oplus\Omega^{1}_{0}(M, {\mathcal P}).
\end{equation}

Next, we show that $\mathrm{Im}(P)$ decomposes as
\begin{equation}\label{cheie}
\mathrm{Im}(P) =  \delta \Omega^{1}_{0}(M, {\mathcal P}) \oplus
C(N({\mathcal P})).
\end{equation}
With the previous notations, it can be checked that $\delta (A) =
P\delta \left(\sum_{i=1}^{3}\alpha_{i}\otimes J_{i}\right) .$ This
implies that $\delta \Omega^{1}_{0}(M, {\mathcal P})$ is included
in $\mathrm{Im}(P)$. We now show that $\delta \Omega^{1}_{0}(M,
{\mathcal P})$ and $C(N({\mathcal P}))$ have trivial intersection.
For this, we need the following observation: for any $B=
\sum_{i=1}^{3}\beta_{i}\otimes J_{i}\in \Omega^{1}_{0}(M,
{\mathcal P})$,
\begin{equation}\label{ul}
\beta_{i}(X)= \frac{\epsilon_{i}\mathrm{tr}\left( J_{i}(\delta
B)_{X}\right)}{n-2},\quad \forall X, \quad\forall i\in \{1 ,2,3\}.
\end{equation}
Relation (\ref{ul}) can be checked using definitions: we first
write the coefficients $\beta_{i}$ of $B$ in the form
\begin{align*}
\beta_{1}&=\frac{2}{3}\gamma_{1}-\frac{1}{3}\gamma_{2}\circ
J_{3}+\frac{1}{3}\gamma_{3}\circ J_{2}\\
\beta_{2}&=\frac{2}{3}\gamma_{2}-\frac{1}{3}\gamma_{1}\circ
J_{3}-\frac{1}{3}\gamma_{3}\circ J_{1}\\
\beta_{3}&= \frac{2}{3}\gamma_{3} +\frac{1}{3}\gamma_{1}\circ
J_{2}+\frac{1}{3}\gamma_{2}\circ J_{1}
\end{align*}
for some $1$-forms $\gamma_{1}$, $\gamma_{2}$, $\gamma_{3}$ and
then we apply the definition of the torsion map $\delta$, we take
traces, etc, and we get (\ref{ul}).  Suppose now that $\delta
(B)\in C(N({\mathcal P}))$. Then, from (\ref{ul}) and the
definition of $C(N({\mathcal P}))$,
$$
0=\mathrm{tr}\left( J_{i} (\delta B)_{X}\right) =
(n-2)\epsilon_{i} \beta_{i}(X),\quad\forall X,\quad \forall i\in
\{ 1,2,3\} ,
$$
which implies that $B=0.$ We proved that
$\delta\Omega^{1}_{0}(M,{\mathcal P})$ and $C(N({\mathcal P}))$
intersect trivially. We now prove that
$\delta\Omega^{1}_{0}(M,{\mathcal P})$ and $C(N({\mathcal P}))$
generate $\mathrm{Im}(P).$ For this, let $T\in \Omega^{2}(M, TM)$
and write
\begin{equation}\label{dt}
P(T)=( P(T)- \delta  ( \sum_{i=1}^{3} \tau_{i}\otimes J_{i})) +
\delta ( \sum_{i=1}^{3} \tau_{i}\otimes J_{i})
\end{equation}
with $1$-forms $\tau_{i}$ defined in (\ref{taui}).  Lemma
\ref{ajut} and the definition of $\tau_{i}$ imply that
$\sum_{i=1}^{3}\tau_{i}\otimes J_{i}$ belongs to
$\Omega^{1}_{0}(M, {\mathcal P})$. Moreover, the first term of
(\ref{dt}) belongs to $C(N({\mathcal P}))$: it belongs to
$\mathrm{Im}(P)$ since $\sum_{i=1}^{3}\tau_{i}\otimes J_{i}\in
\Omega^{1}_{0}(M, {\mathcal P})$ and $\delta\Omega^{1}_{0}(M,
{\mathcal P})$ is included in $\mathrm{Im}(P)$ (from what we
already proved); moreover, for any $i\in \{ 1,2,3\}$ and $X\in
TM$,
$$
\mathrm{tr}\left( J_{i} P(T)_{X}\right)
-\mathrm{tr}\left(J_{i}\delta (\sum_{k=1}^{3}\tau_{k}\otimes
J_{k})_{X}\right) = \mathrm{tr}\left(J_{i}P(T)_{X}\right) -
(n-2)\epsilon_{i}\tau_{i}(X) =0
$$
where the first equality holds from (\ref{ul}), since
$\sum_{i=1}^{3}\tau_{i}\otimes J_{i}\in\Omega^{1}_{0}(M, {\mathcal
P})$, and the second equality is just the definition of
$\tau_{i}.$ The decomposition (\ref{cheie}) follows.

We can now prove the decomposition (\ref{dnq}). Using
(\ref{cheie}) and Lemma \ref{h1}, we obtain the following
decomposition of $\Omega^{2}(M, TM)$:
\begin{equation}\label{o}
\Omega^{2}(M, TM) = \delta\Omega^{1}(M, Z({\mathcal P}))\oplus
\delta \Omega^{1}_{0}(M, {\mathcal P})\oplus C(N({\mathcal P})).
\end{equation}
On the other hand, we claim that
\begin{equation}\label{inn}
\delta\Omega^{1}(M, N({\mathcal P}))= \delta \Omega^{1}(M,
Z({\mathcal P}))\oplus\delta\Omega^{1}_{0}(M, {\mathcal P}),
\end{equation}
or, equivalently, that $\delta \Omega^{1}(M, Z({\mathcal P}))$ and
$\delta\Omega^{1}_{0}(M, {\mathcal P})$ generate
$\delta\Omega^{1}(M, N({\mathcal P}))$ (since $\delta
\Omega^{1}(M, Z({\mathcal P}))$ and $\delta\Omega^{1}_{0}(M,
{\mathcal P})$ have trivial intersection, from (\ref{o})). Recall
that $N({\mathcal P}) = Z({\mathcal P})\oplus {\mathcal P}.$ From
(\ref{initial}), we notice that in order to prove (\ref{inn}) it
is enough to show that $\delta {\mathcal C}(M)$ is included in
$\delta\Omega^{1}(M, Z({\mathcal P}))$. Let $T^{\alpha}=
\sum_{i=1}^{3}\epsilon_{i} (\alpha\circ J_{i})\otimes
J_{i}\in{\mathcal C}(M)$, where $\alpha \in\Omega^{1}(M)$ is an
arbitrary $1$-form. It can be checked that $\delta (T^{\alpha})=
\delta (E^{\alpha})$, where $E^{\alpha}$, defined by
$$
E^{\alpha}(X, Y) :=-\left(  \alpha (Y)X
+\sum_{i=1}^{3}\epsilon_{i} \alpha (J_{i}Y)J_{i}X +\alpha
(X)Y\right) ,\quad\forall X, Y\in TM,
$$
belongs to $\Omega^{1}(M, Z({\mathcal P}))$. This implies that
$\delta {\mathcal C}(M)$ is included in $\delta \Omega^{1}(M,
Z({\mathcal P}))$ as claimed. Decomposition (\ref{inn}) follows
and implies, together with (\ref{o}), decomposition (\ref{dnq}).
Clearly, the map (\ref{proj}) is the projection onto the second
factor of (\ref{dnq}).

\end{proof}

As a consequence, we recover Proposition 2.5 of \cite{zam}.

\begin{cor}
The torsion tensor $T^{\mathcal P}$ of an almost para-quaternionic
structure  $\mathcal P$ is related to the torsion tensor $T^{H}$
of a compatible almost para-hypercomplex structure $H= \{ J_{1},
J_{2}, J_{3}\}$ by
$$
T^{\mathcal P}:=T^{H}-\delta\left( \tau_{1}\otimes
J_{1}+\tau_{2}\otimes J_{2}+\tau_{3}\otimes J_{3}\right) ,
$$
where, for any tangent vector $X$,
\begin{equation}\label{th}
\tau_{i}(X):=\frac{\epsilon_{i}\mathrm{tr}(J_{i}T^{H}_{X})}{n-2},\quad
\forall i\in \{ 1,2,3\}.
\end{equation}
\end{cor}

\begin{proof} The torsion tensor
$T^{\mathcal P}\in C(N({\mathcal P}))$ is the projection of
$T^{H}\in\mathrm{Im}(P)$ with respect to the decomposition
(\ref{cheie}).\end{proof}

We will need the following Corollary in the proof of Theorem
\ref{sug}. This Corollary is analogue to Proposition 2.3 of
\cite{ponte} and can be proved in a similar way. A similar result
has been proved in \cite{zam}.

\begin{cor}\label{fol} The torsion $T^{\mathcal P}$  of an almost para-quaternionic
manifold $(M, {\mathcal P})$ is zero if and only if the torsion
$T^{H}$ of any compatible almost para-hypercomplex structure $H =
\{ J_{1}, J_{2}, J_{3}\}$ is of the form
\begin{equation}\label{inplus}
T^{H} = \delta ( \sum_{i=1}^{3}\alpha_{i}\otimes
J_{i}+\alpha\otimes \mathrm{Id})
\end{equation}
where $\alpha ,\alpha_{1}, \alpha_{2}, \alpha_{3}$ are $1$-forms.
\end{cor}

\subsection{Minimal para-quaternionic connections}

Let $(M, {\mathcal P})$ be an almost para-quaternionic manifold. A
para-quaternionic connection $\nabla$ is minimal if its torsion
$T^{\nabla}$ is equal to the torsion tensor $T^{\mathcal P}\in
C(N({\mathcal P}))$ of the almost para-quaternionic structure
$\mathcal P.$ Minimal para-quaternionic connections always exist
(see \cite{zam}, Proposition 2.5). Moreover, they are parametrized
by $1$-forms, as stated in the following Lemma.

\begin{lem}\label{at}
Let $\{ J_{1}, J_{2}, J_{3}\}$ be an admissible basis of $\mathcal
P .$ Any two minimal para-quaternionic connections $\nabla$ and
$\nabla^{\prime}$ on $(M, {\mathcal P})$ are related by
$\nabla^{\prime}=\nabla+S^{\alpha}$, where $\alpha\in
\Omega^{1}(M)$ and
\begin{align*}
S^{\alpha}_{X}(Y)&:=\alpha
(Y)X-\alpha (J_{1}Y)J_{1}X+\alpha (J_{2}Y)J_{2}X+\alpha (J_{3}Y)J_{3}X\\
&+\alpha (X)Y-\alpha (J_{1}X)J_{1}Y+\alpha (J_{2}X)J_{2}Y+\alpha
(J_{3}X)J_{3}Y,
\end{align*}
for any $X, Y\in TM.$
\end{lem}

\begin{proof}
The statement is more general: two para-quaternionic connections
$\nabla$ and $\nabla^{\prime}$ have the same torsion if and only
if there is a $1$-form $\alpha\in \Omega^{1}(M)$ such that
$\nabla^{\prime} = \nabla + S^{\alpha}$. This comes from the fact
that the map which associates to a covector $\alpha\in T^{*}M$ the
tensor $S^{\alpha}\in T^{*}M\otimes \mathrm{End}(TM)$ defined as
above is an isomorphism between $T^{*}M$ and $\left( T^{*}M
\otimes N({\mathcal P})\right) \cap \left( S^{2}T^{*}M\otimes
TM\right)$, where $S^{2}T^{*}M$ is the bundle of symmetric
$(2,0)$-tensors on $M$.

\end{proof}

\section{Compatible (para)-complex structures}\label{nou}

In our conventions, a system $\{ I_{i}\}$ of almost complex and/or
almost para-complex structures on a manifold $M$ is independent if
it is pointwise independent, i.e. for any $p\in M$, the system $\{
I_{i}(p)\}$ is independent. In particular, two almost complex or
almost para-complex structures $I_{1}$, $I_{2}$ on $M$ are
independent if $I_{1}(p)\neq \pm I_{2}(p)$ at any $p\in M.$

Our main result in this Section is the following criteria of
integrability of almost para-quaternionic structures. Similar
results are known for conformal $4$-manifolds and for almost
quaternionic manifolds (see \cite{ponte}, \cite{ponte1} and
\cite{sal}).

\begin{thm}\label{sug}
Let $(M, {\mathcal P})$ be an almost para-quaternionic manifold of
dimension $4m\geq 8.$ Suppose one of the following situations
holds:
\begin{enumerate}

\item there are two
globally defined, independent, compatible, complex or para-complex
structures $I_{1}$ and $I_{2}$ such that either $I_{1}$ or $I_{2}$
is a complex structure.

\item there are two globally defined, independent, compatible,
para-complex structures $I_{1}$ and $I_{2}$ such that at any $p\in
M$,  the $2$-plane $\mathrm{Span}\{ I_{1}(p), I_{2}(p)\}$ is
non-degenerate  (with respect to the standard Lorenzian metric
$\langle\cdot ,\cdot\rangle$ of $\mathcal P$).

\item there are three globally defined, pairwise independent,
compatible, para-complex structures $I_{1}$, $I_{2}$ and $I_{3}$,
such that at any $p\in M$, $I_{1}(p)$, $I_{2}(p)$ and $I_{3}(p)$
are linearly dependent and for any $i\neq j$, the $2$-plane
$\mathrm{Span}\{ I_{i}(p), I_{j}(p)\}$ is degenerate.

\end{enumerate}

Then $(M, {\mathcal P})$ is para-quaternionic.
\end{thm}

We divide the proof of Theorem \ref{sug} into several Lemmas and
Propositions. We begin with Lemmas \ref{j12} and \ref{nijab},
which are mild generalizations of (3.4.1) and (3.4.4) of
\cite{alek} and can be proved in a similar way.

\begin{lem}\label{j12}
Let $J_{i}$ ($i\in \{ 1,2\}$) be two anti-commuting almost
$\epsilon_{i}$-complex structures on a manifold $M$. The Nijenuis
tensors of $J_{1}$, $J_{2}$ and $J_{1}\circ J_{2}$ are related by
\begin{align*}
2N_{J_{1}\circ J_{2}}(X, Y)= & N_{J_{1}}(J_{2}X, J_{2}Y)
-\epsilon_{1} N_{J_{2}} (X, Y) -J_{1} N_{J_{2}}(J_{1}X, Y)\\
& - J_{1} N_{J_{2}} (X, J_{1}Y) -\epsilon_{2} N_{J_{1}}(X, Y)
+N_{J_{2}} (J_{1}X, J_{1}Y)\\
& -J_{2} N_{J_{1}}(X, J_{2}Y)- J_{2} N_{J_{1}}(J_{2}X, Y),
\end{align*}
for any vector fields $X, Y\in {\mathcal X}(M).$ In particular, if
$J_{1}$ and $J_{2}$ are integrable, also $J_{1}\circ J_{2}$ is.
\end{lem}

Recall that if $A, B\in \mathrm{End}(TM)$ are two endomorphisms of
$TM$, their Nijenuis bracket is a new endomorphism of $TM$,
defined, for any vector fields $X, Y\in {\mathcal X}(M)$, by
\begin{align*}
[A, B](X, Y)&= [AX, BY] +[BX, AY] -A( [BX, Y]+[X, BY])\\
&- B([X,AY]+[AX, Y])+(AB+BA)[X, Y].
\end{align*}
Note that if $J$ is an almost $\epsilon$-complex structure, then
$[J, J] = 2N_{J}.$

\begin{lem}\label{nijab} Let $\nabla$ be a linear connection on a
manifold $M$, which preserves two endomorphisms $A, B\in
\mathrm{End}(TM)$. Denote by $T^{\nabla}$ the torsion of the
connection  $\nabla .$ For any vector fields $X, Y\in {\mathcal
X}(M)$,
\begin{align*}\label{compli}
[A, B](X, Y) &= A\left( T^{\nabla}(BX, Y) + T^{\nabla}(X, BY)\right)\\
&+B\left( T^{\nabla}(X, AY) +T^{\nabla}(AX, Y)\right)\\
&-T^{\nabla}(AX, BY) - T^{\nabla}(BX, AY).
\end{align*}
\end{lem}

In the next Lemma we collect some simple algebraic properties of
compatible almost para-complex structures on  almost
para-quaternionic manifolds.

\begin{lem}\label{baze} Let $(M, {\mathcal P})$ be an almost para-quaternionic
manifold.

i) suppose that $I_{1}$ and $I_{2}$ are two globally defined,
compatible, independent, almost para-complex structures and let
$p\in M$. The $2$-plane $\mathrm{Span}\{ I_{1}(p),I_{2}(p)\}$ is
degenerate if and only if $|\langle I_{1}(p),I_{2}(p)\rangle | =
1$, if and only if $\mathrm{pr}_{I_{1}^{\perp}}(I_{2}) = I_{2}-
\langle I_{2}, I_{1}\rangle I_{1}$ or
$\mathrm{pr}_{I_{2}^{\perp}}(I_{1})=I_{1} -\langle I_{1},
I_{2}\rangle I_{2}$ is null (i.e. squares to the trivial
endomorphism) at the point $p$.

ii) suppose that $I_{1}$, $I_{2}$ and $I_{3}$ are three pairwise
independent, globally defined, compatible, almost para-complex
structures, such that $\mathrm{Span}\{ I_{i}(p),I_{j}(p)\}$ is
degenerate, for any $i\neq j$ and any $p\in M.$ Moreover, assume
that
$$
\langle I_{1}, I_{2}\rangle\langle I_{2}, I_{3}\rangle\langle
I_{1}, I_{3}\rangle =1.
$$
Then the system $\{ I_{1}, I_{2}, I_{3}\}$ is linearly dependent
at any point and (eventually changing the order of $\{ I_{i}\}$
and replacing $I_{i}$ with $-I_{i}$ if necessary) there is, in a
neighborhood of any point, an admissible basis $\{ J_{1}, J_{2},
J_{3}\}$ of $\mathcal P$ such that
\begin{align*}
I_{1} &= J_{2}\\
I_{2}&= J_{1}-J_{2} + qJ_{3}\\
I_{3}&= aJ_{1}+J_{2} +aqJ_{3},
\end{align*}
where $a$ is a smooth function, non-vanishing and different from
$-1$ at any point, and $q\in \{ -1, 1\} .$

iii) suppose that $I_{1}$, $I_{2}$ and $I_{3}$ are like in ii),
but
$$
\langle I_{1}, I_{2}\rangle\langle I_{2}, I_{3}\rangle\langle
I_{1}, I_{3}\rangle =-1.
$$
Then the system $\{ I_{1}, I_{2}, I_{3}\}$ is linearly independent
and (eventually changing the order of $\{I_{i}\}$ and replacing
$I_{i}$ with $-I_{i}$ if necessary) there is a global admissible
basis $\{ J_{1}, J_{2}, J_{3}\}$ of $\mathcal P$ such that
\begin{align*}
I_{1}& = J_{2}\\
I_{2}&= J_{1}+ J_{2}+ J_{3}\\
I_{3}&= J_{1}+ J_{2}-J_{3}.
\end{align*}
\end{lem}

\begin{proof}
The first statement is easy. To prove ii), suppose that $\langle
I_{1}, I_{3}\rangle =1.$ Replacing $I_{2}$ with $-I_{2}$ if
necessary, we can moreover assume that both $\langle I_{1},
I_{2}\rangle$ and $\langle I_{2}, I_{3}\rangle$ are equal to $-1.$
Then $I_{1}+I_{2}$ and $I_{3} -I_{1}$ belong to $I_{1}^{\perp}$,
are null and orthogonal. Therefore, they must be proportional (the
restriction of $\langle\cdot ,\cdot\rangle$ to $I_{1}^{\perp}$
being non-degenerate). We deduce that $\{I_{1}, I_{2}, I_{3}\}$
are dependent at any point. Clearly, we can find an admissible
basis of $\mathcal P$ such that $I_{1}=J_{2}$ and $I_{2}=
J_{1}-J_{2}+qJ_{3}$, where $q\in \{ -1, 1\}.$ Since $I_{1}+I_{2}$
and $I_{3} -I_{1}$ are proportional, $I_{3}=I_{1}+
a(J_{1}+qJ_{3})$ for a smooth function $a$. Since $I_{1}$, $I_{2}$
and $I_{3}$ are pairwise independent, $a$ is non-vanishing and
different from $-1$ at any point. The second claim follows. The
third claim is equally easily.
\end{proof}

Using the previous Lemmas, we can now prove Theorem \ref{sug}. We
first assume in Proposition \ref{integr} that $(M, {\mathcal P})$
admits a pair of (integrable) complex or para-complex structures
like in the first two cases of Theorem \ref{sug} and we show that
$\mathcal P$ is para-quaternionic.  The remaining case of Theorem
\ref{sug} will be treated in Proposition \ref{main1}.

\begin{prop}\label{integr} Let $(M, {\mathcal P})$ be an almost para-quaternionic
manifold of dimension $n=4m\geq 8.$ Suppose that $\mathcal P$
admits two globally defined, independent, compatible, complex or
para-complex structures $I_{1}$ and $I_{2}$, such that one of the
following conditions holds:

i) either $I_{1}$ or $I_{2}$ is a complex structure.

ii) both $I_{1}$ and $I_{2}$ are para-complex structures and
$\mathrm{Span}\{ I_{1}(p), I_{2}(p)\}$ is non-degenerate at any
$p\in M.$

Then $(M, {\mathcal P})$ is a para-quaternionic manifold.
\end{prop}

\begin{proof} Our argument is similar to the one
employed in the proof of Theorem 2.4 of \cite{ponte}. In a
neighborhood of any point we consider two almost
$\epsilon_{i}$-complex structures $J_{i}$ ($i\in \{ 1,2\}$), with
$\langle J_{1}, J_{2}\rangle =0$, such that $I_{1} = J_{1}$ and
$I_{2} = aJ_{1}+bJ_{2}$, where $a, b$ are smooth functions, with
$b$ non-vanishing (this is possible since
$\mathrm{pr}_{I_{1}^{\perp}}(I_{2})$ is non-null when $I_{1}$ is
complex --- the metric on $I_{1}^{\perp}$ being positive definite
--- and also when both $I_{1}$ and $I_{2}$ are para-complex, from the
non-degeneracy of the $2$-planes $\mathrm{Span}\{ I_{1}(p),
I_{2}(p)\}$ and  Lemma \ref{baze}). Since $\langle J_{1},
J_{2}\rangle =0$, $J_{1}$ and $J_{2}$ anti-commute and the
composition $J_{3}:= J_{1}\circ J_{2}$ is an almost
$\epsilon_{3}$-complex structure, with $\epsilon_{3} := -
\epsilon_{1}\epsilon_{2}\in \{ -1,+1\}$. We divide the proof into
three steps.

Step one: we prove that the torsion $T^{H}$ of the almost
para-hypercomplex structure $H$ defined by $J_{1}$, $J_{2}$ and
$J_{3}$ has the following expression: for any vector fields $X$,
$Y$,
\begin{eqnarray*}\label{t}
-12 T^{H}(X, Y) &= &3\epsilon_{2} N_{J_{2}} (X, Y)
-\epsilon_{3}J_{1}N_{J_{2}}(J_{1}X,Y)
-\epsilon_{3}J_{1}N_{J_{2}}(X, J_{1}Y) \\
&&+\epsilon_{3} N_{J_{2}} (J_{1}X, J_{1}Y).
\end{eqnarray*}
We prove this in the following way: since $J_{1}=I_{1}$ is
integrable, Lemma \ref{torsionc} implies that
\begin{equation}\label{need}
T^{H} = -\frac{1}{6}\left( \epsilon_{2}N_{J_{2}}+\epsilon_{3}
N_{J_{3}}\right) .
\end{equation}
From Lemma \ref{j12} and the integrability of $J_{1}$,
\begin{equation}\label{nj123}
2N_{J_{3}}(X, Y)= -\epsilon_{1} N_{J_{2}} (X, Y) -J_{1}
N_{J_{2}}(J_{1}X, Y) - J_{1} N_{J_{2}} (X, J_{1}Y)+N_{J_{2}}
(J_{1}X, J_{1}Y).
\end{equation}
Combining (\ref{need}) with (\ref{nj123}) we get our first claim.

Step two: we prove that the Nijenhuis bracket $[J_{1}, J_{2}]$ has
the following expression:
\begin{align*} 2[ J_{1},J_{2}] (X, Y)= &\epsilon_{2}J_{3} N_{J_{2}}(X, Y) -\epsilon_{3}
J_{3} N_{J_{2}}(J_{1}X, J_{1}Y)\\
& -\epsilon_{2} J_{2} \left( N_{J_{2}}(J_{1}X, Y) +N_{J_{2}}(X,
J_{1}Y)\right) .
\end{align*}
To prove this claim,  we apply Lemma \ref{nijab} to $A:= J_{1}$,
$B:= J_{2}$ and the Obata connection $\nabla$ of $H$ (which
preserves $J_{1}$, $J_{2}$ and $J_{3}$). Since $T^{H}=T^{\nabla}$,
from Lemma \ref{nijab},
\begin{align*}
[J_{1}, J_{2}](X, Y)&= J_{1}\left(  T^{H}(J_{2}X, Y) +T^{H} (X,
J_{2}Y)\right)\\
&+J_{2}\left( T^{H} (J_{1}X,Y)+T^{H}(X, J_{1}Y)\right)\\
& - T^{H}(J_{1}X, J_{2}Y) -T^{H}(J_{2}X, J_{1}Y) .
\end{align*}
We now evaluate $T^{H}(J_{2}X, Y) +T^{H} (X, J_{2}Y)$, $T^{H}
(J_{1}X, Y)+T^{H}(X, J_{1}Y)$ and $T^{H}(J_{1}X, J_{2}Y)
+T^{H}(J_{2}X, J_{1}Y)$ in terms of the Nijenhuis tensor
$N_{J_{2}}.$ It can be checked that the Nijenuis tensor of any
almost $\epsilon$-complex structure $J$ has the following
symmetries:
\begin{equation}\label{n}
N_{J} (JX,Y) = N_{J}(X, JY) = - JN_{J} (X, Y),\quad\forall X, Y\in
{\mathcal X}(M).
\end{equation}
Using the expression of $T^{H}$ determined in the first step,
relation (\ref{n}) for $J := J_{2}$ and the anti-commutativity
$J_{1}J_{2}= - J_{1}J_{2}$ we obtain
\begin{align*}
6\left( T^{H}(J_{1}X,Y) + T^{H}(X, J_{1}Y)\right)  =& \epsilon_{3}
J_{1} N_{J_{2}} (J_{1}X, J_{1}Y) -\epsilon_{2}
J_{1}N_{J_{2}} (X, Y)\\
& -\epsilon_{2} \left( N_{J_{2}}(J_{1}X, Y) +N_{J_{2}}(X,
J_{1}Y)\right)\\
6\left( T^{H}(J_{2}X, Y) +T^{H}(X, J_{2}Y)\right)=& 3\epsilon_{2}
J_{2} N_{J_{2}} (X, Y) -\epsilon_{3} J_{2}
N_{J_{2}} (J_{1}X, J_{1}Y)\\
6\left( T^{H}(J_{1}X, J_{2}Y) + T^{H}(J_{2}X, J_{1}Y)\right) &=
2\epsilon_{2}J_{2}\left( N_{J_{2}}(J_{1}X, Y) + N_{J_{2}}(X,
J_{1}Y)\right)\\
& +\epsilon_{3}J_{3}\left( N_{J_{2}}(J_{1}X, J_{1}Y) -
\epsilon_{1}N_{J_{2}}(X, Y)\right) .
\end{align*}
Replacing these relations in the expression of $[J_{1}, J_{2}]$
above we get our second claim.

Step three: we prove that $T^{H}\equiv 0$, where the sign
"$\equiv$" means equality, modulo terms of the form $\delta
(\alpha_{1}\otimes J_{1}+\alpha_{2}\otimes J_{2}+\alpha_{3}\otimes
J_{3}+\alpha\otimes \mathrm{Id})$, where $\alpha ,\alpha_{1},
\alpha_{2}, \alpha_{3}$ are $1$-forms. To prove this claim, we
notice that, since $I_{2} = a J_{1}+bJ_{2}$, is integrable
\begin{equation}\label{ni2}
a^{2} N_{J_{1}} + b^{2} N_{J_{2}} +ab [ J_{1}, J_{2}]\equiv 0.
\end{equation}
The integrability of $J_{1}$ together with (\ref{ni2}) imply that
\begin{equation}\label{star}
b^{2} N_{J_{2}} +ab [ J_{1}, J_{2}]\equiv 0.
\end{equation}
On the set of points $M_{0}\subset M$ where $a=0$, (\ref{star})
implies that $N_{J_{2}}\equiv 0$ (because $b$ is non-vanishing).
We use now the expression of $[J_{1}, J_{2}]$ determined in  Step
two to show that (\ref{star}) implies that $N_{J_{2}}\equiv 0$
also on $M\setminus M_{0}.$ On $M\setminus M_{0}$, we can divide
(\ref{star}) by $a$ and, using the expression of $[J_{1}, J_{2}]$
provided by Step two, we obtain
\begin{align*}
-\frac{2b}{a} N_{J_{2}} (X, Y)\equiv & \epsilon_{2} J_{3}
N_{J_{2}}(X,
Y) -\epsilon_{3} J_{3} N_{J_{2}} (J_{1}X, J_{1}Y)\\
&-\epsilon_{2} J_{2}\left( N_{J_{2}}(J_{1}X, Y) +N_{J_{2}}(X,
J_{1}Y)\right) .
\end{align*}
Replacing $(X, Y)$ with $(J_{2}X, J_{2}Y)$ in this relation, using
again (\ref{n}) for $J = J_{2}$ and the anti-commutativity
$J_{1}J_{2}= - J_{2}J_{1}$ we get two relations:
\begin{equation}\label{r1}
-\frac{2b}{a} N_{J_{2}}(X, Y)\equiv \epsilon_{2}J_{3}
N_{J_{2}}(X,Y) -\epsilon_{3}J_{3} N_{J_{2}}(J_{1}X, J_{1}Y)
\end{equation}
and
\begin{equation}\label{r2}
N_{J_{2}}(J_{1}X, Y)+ N_{J_{2}}(X, J_{1}Y)\equiv 0.
\end{equation}
Relation (\ref{r2}) implies that
\begin{equation}\label{clar}
N_{J_{2}}(J_{1}X, J_{1}Y)\equiv -\epsilon_{1} N_{J_{2}}(X, Y).
\end{equation}
Replacing (\ref{clar}) in (\ref{r1}) and using $\epsilon_{3} = -
\epsilon_{1}\epsilon_{2}$ we get $N_{J_{2}}\equiv 0$.  From
(\ref{nj123}), $N_{J_{3}}\equiv 0$ as well and then, from
(\ref{t}), $T^{H}\equiv 0.$ Corollary \ref{fol} implies that
$T^{\mathcal P}=0$. This concludes our proof.

\end{proof}

\begin{prop}\label{main1} Let $(M, {\mathcal P})$ be an almost para-quaternionic
manifold of dimension $n=4m\geq 8.$ Suppose that $\mathcal P$
admits three pairwise independent, compatible, para-complex
structures $\{ I_{1}, I_{2}, I_{3}\}$, such that at any $p\in M$,
$\{ I_{1}(p), I_{2}(p), I_{3}(p)\}$ are dependent and for any
$i\neq j$, $\mathrm{Span}\{ I_{i}(p), I_{j}(p)\}$ is degenerate.
Then $(M, {\mathcal P})$ is para-quaternionic.
\end{prop}

\begin{proof}
Like in the proof of Proposition \ref{integr}, we will determine,
in a neighborhood of any point, a suitable compatible almost
para-hypercomplex structure $H$, for which $T^{H}\equiv 0.$ We
divide the proof into two steps.

Step one: let $H:=\{J_{1}, J_{2}, J_{3}\}$ be any admissible basis
of $\mathcal P$ such that $I_{1}= J_{2}$. In particular, $J_{2}$
is integrable. We claim that the torsion $T^{H}$ of $H$ and the
Nijenhuis brackets $[J_{1}, J_{2}]$, $[J_{1}, J_{3}]$ and $[J_{2},
J_{3}]$ have the following expressions: for any vector fields $X$
and $Y$,
\begin{equation}\label{torsiunea}
12 T^{H}(X, Y)= 3 N_{J_{1}}(X, Y) - N_{J_{1}}(J_{2}X, J_{2}Y)
 + J_{2}\left( N_{J_{1}}(X, J_{2}Y) +N_{J_{1}}(J_{2}X, Y)\right)
\end{equation}
and
\begin{align*}
2[J_{1}, J_{2}](X, Y)&= J_{1}\left( N_{J_{1}}(J_{2}X, Y) +
N_{J_{1}}(X, J_{2}Y)\right)\\
&+J_{3}\left( N_{J_{1}}(J_{2}X, J_{2}Y) + N_{J_{1}}(X, Y)\right)\\
2 [J_{1}, J_{3}](X, Y)&= N_{J_{1}}(J_{2}X, Y) + N_{J_{1}}(X,
J_{2}Y)\\
& -J_{2}\left( N_{J_{1}}(J_{2}X, J_{2}Y)+N_{J_{1}}(X,Y)\right)\\
2[J_{2}, J_{3}](X, Y)&= J_{1}\left( N_{J_{1}}(J_{2}X, J_{2}Y) +
N_{J_{1}}(X, Y)\right)\\
&+ J_{3}\left( N_{J_{1}}(J_{2}X, Y) + N_{J_{1}}(X, J_{2}Y)\right).
\end{align*}
To prove these claims, notice that
\begin{equation}\label{tors}
T^{H}= -\frac{1}{6}\left( -N_{J_{1}}+N_{J_{3}}\right)
\end{equation}
since $J_{2}$ is integrable. On the other hand, applying Lemma
\ref{j12} to $J_{1}$ and $J_{2}$ and using the integrability of
$J_{2}$ we get
\begin{equation}\label{nun}
2N_{J_{3}}(X, Y) = N_{J_{1}}(J_{2}X, J_{2}Y) - N_{J_{1}}(X, Y)
-J_{2} N_{J_{1}}(X, J_{2}Y) - J_{2}N_{J_{1}}(J_{2}X, Y).
\end{equation}
Combining (\ref{tors}) with (\ref{nun}) we obtain
(\ref{torsiunea}). In order to evaluate the Nijenhuis brackets
$[J_{1}, J_{2}]$, $[J_{1}, J_{3}]$ and $[J_{2}, J_{3}]$ we will
use Lemma \ref{nijab}, with $\nabla$ the Obata connection of $\{
J_{1}, J_{2}, J_{3}\}$, so that $T^{\nabla} = T^{H}.$ Using
(\ref{torsiunea}), the anti-commutativity of $J_{1}$, $J_{2}$,
$J_{3}$ and (\ref{n}) for $J = J_{1}$, we obtain
\begin{align*}
6\left( T^{H}(J_{1}X, J_{2}Y) + T^{H}(J_{2}X, J_{1}Y)\right) &=
-2J_{1}\left( N_{J_{1}}(X, J_{2}Y)  + N_{J_{1}}(J_{2}X, Y)\right)\\
&-J_{3}\left( N_{J_{1}}(J_{2}X, J_{2}Y) - N_{J_{1}}(X, Y)\right)\\
6\left( T^{H}(J_{1}X, J_{3}Y)  + T^{H}(J_{3}X, J_{1}Y)\right)& =
-N_{J_{1}}(X, J_{2}Y) - N_{J_{1}}(J_{2}X, Y)\\
&+J_{2}\left( N_{J_{1}}(J_{2}X, J_{2}Y)  + N_{J_{1}}(X,
Y)\right)\\
6\left( T^{H}(J_{2}X, J_{3}Y) + T^{H}(J_{3}X, J_{2}Y)\right)& =
-J_{1}\left( 3 N_{J_{1}}(J_{2}X, J_{2}Y) + N_{J_{1}}(X, Y)\right)\\
6\left(T^{H} (J_{1}X,Y)+ T^{H}(X, J_{1}Y)\right) &= -J_{1}\left(
3N_{J_{1}}(X,Y) + N_{J_{1}}(J_{2}X, J_{2}Y)\right)\\
6\left( T^{H}(J_{2}X, Y)+ T^{H}(X, J_{2}Y)\right) &= J_{2}\left(
N_{J_{1}}(J_{2}X, J_{2}Y) + N_{J_{1}}(X, Y)\right)\\
& +N_{J_{1}}(J_{2}X, Y) + N_{J_{1}}(X, J_{2}Y)\\
6 \left( T^{H}(X, J_{3}Y) + T^{H}(J_{3}X, Y)\right) & =
-2J_{1}\left( N_{J_{1}}(J_{2}X, Y) + N_{J_{1}}(X, J_{2}Y)\right)\\
&+ J_{3}\left( N_{J_{1}}(J_{2}X, J_{2}Y) - N_{J_{1}}(X,
Y)\right) .\\
\end{align*}
Applying  Lemma \ref{nijab} we get the expressions of $[J_{1},
J_{2}]$, $[J_{1}, J_{3}]$ and $[J_{2}, J_{3}]$ as stated.

Step two: using Step one, we prove that $\mathcal P$ is
integrable. From Lemma \ref{baze} we can chose, in a neighborhood
of any point, an admissible basis $\{ J_{1}, J_{2}, J_{3}\}$ of
$\mathcal P$ such that
$$
I_{1}= J_{2},\quad I_{2}= J_{1}-J_{2}+qJ_{3},\quad I_{3}=
aJ_{1}+J_{2}+aqJ_{3},
$$
where $a$ is a smooth function, non-vanishing and different from
$-1$ at any point, and $q\in \{ -1,+1\} .$ The integrability of
$I_{1}$ and $I_{3}$ implies that
\begin{align*}
N_{J_{1}}+ N_{J_{3}}+q [J_{1}, J_{3}]&\equiv 0\\
[J_{1}, J_{2}]+q[J_{2}, J_{3}]&\equiv 0,
\end{align*}
where the sign $\equiv$ has the same meaning as in the proof of
Proposition \ref{integr}. Using (\ref{nun}) and the expressions of
$[J_{1}, J_{2}]$ and $[J_{2}, J_{3}]$ previously determined, we
get
\begin{align*}
(\mathrm{Id}&+q J_{2})\{ N_{J_{1}}(J_{2}X, J_{2}Y) + N_{J_{1}}(X,
Y)) + q N_{J_{1}}(J_{2}X, Y)+ q N_{J_{1}}(X, J_{2}Y)\}
\equiv 0\\
(\mathrm{Id} &-q J_{2})\{ N_{J_{1}}(J_{2}X, J_{2}Y) + N_{J_{1}}(X,
Y)) + qN_{J_{1}}(J_{2}X, Y)+ q N_{J_{1}}(X,
J_{2}Y)\} \equiv 0.\\
\end{align*}
Adding these relations,
\begin{equation}\label{adaugat}
N_{J_{1}}(J_{2}X, J_{2}Y) + N_{J_{1}}(X, Y) + q(N_{J_{1}}(J_{2}X,
Y)+ N_{J_{1}}(X, J_{2}Y))\equiv 0.
\end{equation}
Replacing in (\ref{adaugat}) the pair $(X, Y)$ with $(J_{1}X,
J_{1}Y)$ and using (\ref{n}) for $J = J_{1}$ we get
\begin{align*}
N_{J_{1}}(J_{2}X, J_{2}Y) &+ N_{J_{1}}(X, Y)\equiv 0\\
N_{J_{1}}(J_{2}X, Y)&+ N_{J_{1}}(X, J_{2}Y)\equiv 0.
\end{align*}
Like in the proof of Proposition \ref{integr}, we deduce that
$N_{J_{1}}\equiv 0.$ From (\ref{nun}) it follows that
$N_{J_{3}}\equiv 0$ as well and therefore, $T^{H}\equiv 0.$ We
conclude that $\mathcal P$ is para-quaternionic.

\end{proof}

Proposition \ref{main1} concludes the proof of Theorem
\ref{sug}.\\

Theorem \ref{sug} raises the following question: does the
existence of three globally defined, independent, compatible,
para-complex structures $\{ I_{1}, I_{2}, I_{3}\}$ on an almost
para-quaternionic manifold $(M, {\mathcal P})$, such that for any
$p\in M$ and $i\neq j$, the $2$-plane $\mathrm{Span}\{ I_{i}(p),
I_{j}(p)\}$ is degenerate, imply the integrability of the almost
para-quaternionic structure $\mathcal P$? We will now show that
the answer to this question is negative.

For this, it is convenient to express the integrability of
$I_{1}$, $I_{2}$ and $I_{3}$ in terms of the admissible basis
$\{J_{1}, J_{2}, J_{3}\}$ of $\mathcal P$ provided by Lemma
\ref{baze}, i.e. related to $\{ I_{1}, I_{2}, I_{3}\}$ by
\begin{equation}\label{i1i2i3}
I_{1}=J_{2},\quad I_{2}= J_{1}+J_{2}+J_{3},\quad I_{3}=
J_{1}+J_{2}-J_{3}.
\end{equation}

\begin{lem}\label{pe} The integrability of the almost para-complex
structures $I_{1}$, $I_{2}$ and $I_{3}$ is equivalent to the
integrability of $J_{2}$ together with the integrability of the
eigenbundle of $J_{3}$ which corresponds to the eigenvalue $+1.$
\end{lem}

\begin{proof} From (\ref{i1i2i3}),  the integrability of $I_{1}$,
$I_{2}$ and $I_{3}$ is equivalent to the integrability of $J_{2}$
and the following two relations:
\begin{align*}
N_{J_{1}}+N_{J_{3}} + [J_{1}, J_{2}]&=0\\
[J_{1}, J_{3}] + [J_{2}, J_{3}]&=0.
\end{align*}
We now express $N_{J_{3}}$, $[J_{1}, J_{2}]$, $[J_{1}, J_{3}]$ and
$[J_{2}, J_{3}]$ in terms of $N_{J_{1}}$, using the computations
of the proof of Proposition \ref{main1}. To simplify notations,
define, for any vector fields $X$ and $Y$,
\begin{equation}\label{e}
E(X, Y):= N_{J_{1}}(J_{2}X, J_{2}Y) + N_{J_{1}}(X, Y)
\end{equation}
and
\begin{equation}\label{f}
F(X, Y) :=N_{J_{1}}(J_{2}X, Y) + N_{J_{1}}(X, J_{2}Y).
\end{equation}
An easy argument shows that the previous two relations reduce to
\begin{equation}\label{single}
(\mathrm{Id}+ J_{3})E(X, Y)=0\quad\forall X, Y.
\end{equation}
Multiplying (\ref{single}) by $J_{2}$ on the left and using
(\ref{n}) for $J = J_{1}$ we obtain
\begin{equation}\label{mai}
J_{2} \left( N_{J_{1}}(J_{2}X, J_{2}Y) + N_{J_{1}}(X, Y)\right) =
N_{J_{1}}(J_{2}J_{1}X, J_{2}Y)- N_{J_{1}}(J_{1}X, Y).
\end{equation}
Replacing in (\ref{mai}) $X$ with $J_{1}X$ we obtain
\begin{equation}\label{mai2}
J_{2} \left( N_{J_{1}}(J_{2}J_{1}X, J_{2}Y) + N_{J_{1}}(J_{1}X,
Y)\right) = - N_{J_{1}}(J_{2}X, J_{2}Y)+ N_{J_{1}}(X, Y).
\end{equation}
On the other hand, multiplying (\ref{mai}) on the left with
$J_{2}$ we obtain
\begin{equation}\label{mai3}
J_{2} (N_{J_{1}}(J_{2}J_{1}X, J_{2}Y) - N_{J_{1}}(J_{1}X, Y)) =
N_{J_{1}}(J_{2}X, J_{2}Y) + N_{J_{1}}(X, Y).
\end{equation}
Using (\ref{mai2}) and (\ref{mai3}), we get that (\ref{single}) is
equivalent to
\begin{equation}\label{nj1j2j33}
N_{J_{1}}(J_{3}X, J_{3}Y)= J_{3}N_{J_{1}}(X, Y),\quad\forall X, Y.
\end{equation}
Using (\ref{nj1j2j33}), we easily get our claim: relation
(\ref{nun}) expresses $N_{J_{3}}$ in terms of $N_{J_{1}}$;
conversely, using Lemma \ref{j12} and the integrability of
$J_{2}$, we can express $N_{J_{1}}$ in terms of $N_{J_{3}}$ as
\begin{equation}\label{new}
2N_{J_{1}}(X, Y)= - N_{J_{3}}(X, Y) - J_{2} N_{J_{3}}(J_{2}X, Y) -
J_{2}N_{J_{3}}(X, J_{2}Y) + N_{J_{3}}(J_{2}, J_{2}Y).
\end{equation}
Using (\ref{nun}) and (\ref{new}), it can be checked that
(\ref{nj1j2j33}) is equivalent to
$$
J_{3}N_{J_{3}}(X, Y)= N_{J_{3}}(X, Y),\quad\forall X, Y,
$$
or to the integrability of the eigenbundle of $J_{3}$ which
corresponds to the eigenvalue $+1$. Our claim follows.
\end{proof}

Lemma \ref{pe} reduces the problem of finding independent,
compatible, para-complex structures $I_{1}$, $I_{2}$ and $I_{3}$
on $(M, {\mathcal P})$, with $\mathrm{Span}\{ I_{i}(p),
I_{j}(p)\}$ degenerate at any $p\in M$, for any $i\neq j$, to the
problem of finding admissible bases $\{ J_{1}, J_{2}, J_{3}\}$ of
$\mathcal P$, such that both $J_{2}$ and $\mathrm{Ker}
(J_{3}-\mathrm{Id})$ are integrable. We now show that these
conditions, in turn, reduce to solving a certain system of partial
differential equations. Indeed,  since $J_{2}$ is integrable,
locally $M = M_{1}\times M_{2}$ is a product manifold and
$TM_{1}$, $TM_{2}$ are the eigenbundles of $J_{2}$ which
correspond to the eigenvalues $1$ and $-1$, respectively. The four
distributions ${\mathcal
D}^{\pm}:=\mathrm{Ker}(J_{3}\mp\mathrm{Id})$, $TM_{1}$ and
$TM_{2}$ are pairwise transverse. Therefore,
$$
{\mathcal D}^{\pm}:= \{ X+A^{\pm}X,\quad\forall X\in TM_{1}\}
$$
where $A^{+},A^{-}, A^{+}- A^{-}: TM_{1}\rightarrow TM_{2}$ are
isomorphisms. It is easy to check that $J_{2}\circ J_{3} =
-J_{3}\circ J_{2}$ is equivalent to $A^{+}+ A^{-}=0.$ Note also
that in terms of $A:= A^{+}$,
$$
J_{1}(X) = J_{3}(X) = A(X),\quad J_{1}(Y) = -J_{3}(Y) =
-A^{-1}(Y),
$$
for any $X\in TM_{1}$ and $Y\in TM_{2}.$ Let $(x_{1},\cdots,
x_{2m})$ and $(y_{1}, \cdots , y_{2m})$ be local coordinates on
$M_{1}$ and $M_{2}$ respectively. In these coordinates,
$$
A\left(\frac{\partial}{\partial x_{i}}\right)
:=\sum_{j=1}^{2m}f_{ij}\frac{\partial}{\partial y_{j}},
$$
for some smooth functions $f_{ij}= f_{ij}(x_{s},y_{r})$ ($1\leq
i,j, r,s\leq 2m$), with $\mathrm{det}(f_{ij})\neq 0$ at any point.
The system of partial differential equations mentioned above comes
from the integrability of ${\mathcal
D}^{+}=\mathrm{Ker}(J_{3}-\mathrm{Id})$: it can be checked that
${\mathcal D}^{+}$ is integrable if and only if
\begin{equation}\label{ecua}
\frac{\partial f_{jt}}{\partial x_{i}}- \frac{\partial
f_{it}}{\partial x_{j}} +\sum_{k=1}^{2m} \left(
f_{ik}\frac{\partial f_{jt}}{\partial y_{k}} - f_{jk}
\frac{\partial f_{it}}{\partial y_{k}}\right) =0,\quad \forall
i\neq j,\quad \forall t.
\end{equation}

Reversing this argument, any solution $(f_{ij})$ of (\ref{ecua}),
with  non-vanishing determinant $\mathrm{det}(f_{ij})$, defines an
almost para-quaternionic structure ${\mathcal
P}_{(f_{ij})}:=\mathrm{Span} \{ J_{1}, J_{2}, J_{3}\}$, which
admits three independent, compatible, para-complex structures
$I_{1}$, $I_{2}$ and $I_{3}$, related to $\{ J_{1}, J_{2},
J_{3}\}$ by (\ref{i1i2i3}). The next Proposition determines a
class of solutions of (\ref{ecua}), for which $\mathcal P$ is not
para-quaternionic.

\begin{prop}\label{propo} The functions $f_{ij}:=
f_{i}\delta_{ij}$ ($1\leq i,j\leq 2m$), where
$$
f_{1}:=
h\left(\frac{\sum_{j=2}^{2m}x_{j}^{2}}{\sum_{j=2}^{2m}y_{j}^{2}}\right),
\quad
f_{i}:=\frac{x_{i}\left(\sum_{j=2}^{2m}y_{j}^{2}\right)}{y_{i}
\left(\sum_{j=2}^{2m}x_{j}^{2}\right)},\quad 2\leq i\leq 2m
$$
and $h$ is a smooth real function, is a solution of (\ref{ecua}).
Moreover, on any open connected subset of $\mathbb{R}^{4m}$ on
which $f_{i}$ ($1\leq i\leq 2m$) are non-vanishing and
\begin{equation}\label{condities}
h\left(\frac{\sum_{j=2}^{2m}x_{j}^{2}}{\sum_{j=2}^{2m}y_{j}^{2}}\right)
+ \frac{\sum_{j=2}^{2m}x_{j}^{2}}{\sum_{j=2}^{2m}y_{j}^{2}}
\dot{h}\left(\frac{\sum_{j=2}^{2m}x_{j}^{2}}{\sum_{j=2}^{2m}y_{j}^{2}}\right)
\neq 0,
\end{equation}
the associated almost para-quaternionic structure $\mathcal
P_{(f_{ij})}$ is not para-quaternionic.
\end{prop}

\begin{proof} It is straightforward to check that $(f_{ij})$ is a
solution of (\ref{ecua}). We now show that $\mathcal P$ is not
para-quaternionic. With the previous notations, it can be checked
that, for any $i, j\in \{ 1,\cdots , 2m\}$,
\begin{equation}\label{thp}
3 T^{H}\left(\frac{\partial}{\partial x_{i}},
\frac{\partial}{\partial x_{j}}\right) =
\frac{1}{f_{j}}\frac{\partial f_{j}}{\partial x_{i}}
\frac{\partial}{\partial x_{j}} - \frac{1}{f_{i}} \frac{\partial
f_{i}}{\partial x_{j}}\frac{\partial}{\partial x_{i}} ,
\end{equation}
where $T^{H}$ is the torsion of $H:=\{ J_{1}, J_{2}, J_{3}\}$.
Suppose now, by absurd, that $\mathcal P$ is para-quaternionic.
Then $T^{H}$ is of the form
\begin{equation}\label{g1g2g3}
T^{H}=\delta (\gamma\otimes \mathrm{Id}+ \gamma_{1}\otimes J_{1}+
\gamma_{2}\otimes J_{2}+\gamma_{3}\otimes J_{3})
\end{equation}
for some $1$-forms $\gamma$, $\gamma_{1}$, $\gamma_{2}$ and
$\gamma_{3}.$ Moreover, since
$$
J_{1}\left(\frac{\partial}{\partial x_{i}}\right) = J_{3}
\left(\frac{\partial}{\partial x_{i}}\right)
=f_{i}\frac{\partial}{\partial y_{i}},\quad
J_{1}\left(\frac{\partial}{\partial y_{j}}\right) = -J_{3}
\left(\frac{\partial}{\partial y_{j}}\right) =
-\frac{1}{f_{j}}\frac{\partial}{\partial x_{j}},
$$
for any $1\leq i, j\leq 2m$, relation (\ref{g1g2g3}) implies:
\begin{align*}
T^{H}\left(\frac{\partial}{\partial x_{i}},
\frac{\partial}{\partial x_{j}}\right) &= (\gamma
+\gamma_{2})\left(\frac{\partial}{\partial
x_{i}}\right)\frac{\partial}{\partial x_{j}} +
f_{j}(\gamma_{1}+\gamma_{3})\left(\frac{\partial}{\partial
x_{i}}\right)\frac{\partial}{\partial y_{j}}\\
&-(\gamma +\gamma_{2})\left(\frac{\partial}{\partial
x_{j}}\right)\frac{\partial}{\partial x_{i}}
-f_{i}(\gamma_{1}+\gamma_{3})\left(\frac{\partial}{\partial
x_{j}}\right)\frac{\partial}{\partial y_{i}}.
\end{align*}
From (\ref{thp}) we deduce that $(\gamma_{1}+\gamma_{3})
\left(\frac{\partial}{\partial x_{i}}\right) =0$ for any $i$ and
$$
\frac{1}{f_{i}}\frac{\partial f_{i}}{\partial x_{j}} = -(\gamma
+\gamma_{2})\left(\frac{\partial}{\partial x_{j}}\right)
,\quad\forall i\neq j.
$$
Equivalently, for any $i\neq k$, the quotient
$\frac{f_{i}}{f_{k}}$ depends only on $x_{i}, x_{k}, y_{1},\cdots
, y_{2n}.$ However, from the definition of the functions $f_{i}$,
this is cannot hold: take $i=1$, $k\geq 2$ arbitrary and use
(\ref{condities}). We obtain a contradiction. We deduce that
$\mathcal P$ is not para-quaternionic.

\end{proof}

\section{Twistor and reflector spaces}\label{twistor}

Let $(M, {\mathcal P})$ be an almost para-quaternionic manifold.
Denote by $\pi^{\epsilon}:Z^{\epsilon}\rightarrow M$ the reflector
bundle of $(M, {\mathcal P})$ when $\epsilon =1$ and the twistor
bundle of $(M, {\mathcal P})$ when $\epsilon =-1$. A
para-quaternionic connection $\nabla$ on $(M,{\mathcal P})$
induces an almost $\epsilon$-complex structure ${\mathcal
J}^{\nabla ,\epsilon}$ on $Z^{\epsilon}$ as follows: let
$H^{\nabla , \epsilon }_{J}$ be the horisontal space at $J^{}\in
Z^{\epsilon}$ of the connection $\nabla$, acting on the bundle
$\pi^{\epsilon}: Z^{\epsilon}\rightarrow M$, and
\begin{equation}\label{h}
T_{J}Z^{\epsilon}=H^{\nabla , \epsilon}_{J}\oplus
T^{v}_{J}Z^{\epsilon}
\end{equation}
the induced decomposition of $T_{J}Z^{\epsilon}$ into horizontal
and vertical subspaces. On $H^{\nabla , \epsilon}_{J}$, identified
with $T_{\pi^{\epsilon}(J)}M$ by means of the differential
$(\pi^{\epsilon})_{*}$, ${\mathcal J}^{\nabla ,\epsilon}$
coincides with $J$ (viewed as an endomorphism of
$T_{\pi^{\epsilon}(J)}M$); on $T^{v}_{J}Z^{\epsilon}$, it is
$$
{\mathcal J}^{\epsilon}_{J}(A):= J\circ A,\quad\forall A\in
T_{J}^{v}Z^{\epsilon}
$$
and is well defined, since
$$
T^{v}_{J}Z^{\epsilon}= \{ A\in {\mathcal
P}_{\pi^{\epsilon}(J)}:\quad A\circ J + J\circ A=0\}\subset
\mathrm{End} ( T_{\pi^{\epsilon}(J)}M) .
$$
Consider now another para-quaternionic connection
$\nabla^{\prime}.$ If $\nabla$ and $\nabla^{\prime}$ have the same
torsion, ${\mathcal J}^{\nabla ,\epsilon}={\mathcal J}^{
\nabla^{\prime},\epsilon}$ (see Corollary 3.4 of \cite{ivanov}).

\begin{defn} Let $(M,{\mathcal P})$ be an almost para-quaternionic
manifold. The twistor space $Z^{-}$ has a canonical almost complex
structure ${\mathcal J}^{-}:={\mathcal J}^{\nabla ,-}$. Similarly,
the reflector space $Z^{+}$ has a canonical almost para-complex
structure ${\mathcal J}^{+}:={\mathcal J}^{\nabla , +}$. Here
$\nabla$ is (any) minimal para-quaternionic connection on $(M,
{\mathcal P})$.
\end{defn}

We need the following Lemma in the proof of Theorem \ref{final}.

\begin{lem}\label{preg1} Let $(M, {\mathcal P})$ be a para-quaternionic
manifold and $J$ a compatible almost $\epsilon$-complex structure
on $(M, {\mathcal P}).$ Then $J$ is integrable if and only if the
image $\mathrm{Im}(\sigma^{J})\subset Z^{\epsilon}$ of the
tautological section $\sigma^{J}: M\rightarrow Z^{\epsilon}$
defined by $J$ is ${\mathcal J}^{\epsilon}$-stable.
\end{lem}

\begin{proof} The proof goes like in the quaternionic case,
(see \cite{ponte}, Sections 3 and 4), with the Oproiu connections
of almost quaternionic manifolds replaced by the minimal
para-quaternionic connections of $(M,{\mathcal P}).$
\end{proof}

The next Theorem is our main result of this Section and is the
para-quaternionic analogue of Theorem 4.2 of \cite{ponte}. The
equivalence between the second and the third conditions bellow has
already been proved in  \cite{ivanov}.

\begin{thm}\label{final}
Let $(M, {\mathcal P})$ be an almost para-quaternionic manifold of
dimension $n= 4m\geq 8.$ Denote by ${\mathcal J}^{+}$ the
canonical almost para-complex structure of the reflector space
$Z^{+}$ and by ${\mathcal J}^{-}$ the canonical almost complex
structure of the twistor space $Z^{-}$ of $(M, {\mathcal P}).$ The
following statements are equivalent:

i) $(M, {\mathcal P})$ is a para-quaternionic manifold.

ii) both ${\mathcal J}^{+}$ and ${\mathcal J}^{-}$ are integrable.

iii) either ${\mathcal J}^{-}$ or ${\mathcal J}^{+}$ is
integrable.

iv) for any point $p\in M$ and compatible $\epsilon$-complex
structure $I_{p}\in {\mathcal P}_{p}$, there are infinitely many
compatible $\epsilon$-complex structures defined in a neighborhood
of $p$ which extend $I_{p}.$

v) any point of $M$ has a neighborhood on which there are defined
four compatible, pairwise independent, $\epsilon_{i}$-complex
structures $ I_{i}$ ($i\in \{ 1,2,3,4\}$).
\end{thm}

\begin{proof}
Let $\nabla$ be a minimal para-quaternionic connection on $(M,
{\mathcal P})$ , so that $T^{\nabla} = T^{\mathcal P}$ and
${\mathcal J}^{\epsilon} ={\mathcal J}^{\nabla ,\epsilon}$, for
$\epsilon \in \{ -1,1\} .$

We show the equivalence of the first three conditions. Suppose
that $\mathcal P$ is a para-quaternionic structure. Then
$T^{\mathcal P}=0$, the connection $\nabla$ is torsion free and
both $\mathcal J^{+}$, ${\mathcal J}^{-}$ are integrable (see
\cite{ivanov}, Theorem 3.8). Conversely, suppose that $\mathcal
J^{-}$ or ${\mathcal J}^{+}$ is integrable. Then, again from
Theorem 3.8 of \cite{ivanov},
\begin{equation}\label{t02}
\Pi^{0,2}_{J}(T^{\nabla})=0,\quad \forall J\in {\mathcal P} ,\quad
J^{2}=\pm\mathrm{Id}.
\end{equation}
Relation (\ref{t02}) implies that  $T^{\mathcal
P}=P(T^{\nabla})=0$ (where $P$ is the projector of Lemma
\ref{h1}). It follows that $\mathcal P$ is para-quaternionic.

We now show that the fourth and the fifth conditions are
equivalent to any of the first three conditions. From Theorem
\ref{sug}, the fifth condition implies the first. Suppose now that
the first (hence also the second and third) condition holds. We
will prove the fourth condition when $I_{p}$ is para-complex (the
argument when $I_{p}$ is a complex structure is similar). Since
${\mathcal J}^{+}$ is integrable, the distributions
$T^{\pm}Z^{+}:= \mathrm{Ker} ({\mathcal J}^{+}\mp \mathrm{Id})$
are involutive. Being transversal, there are local coordinates
$(x_{1},\cdots ,x_{n+2})$ in a neighborhood $\mathcal U$ of
$I_{p}\in Z^{+}$ such that
\begin{align*}
T^{+}Z^{+}&=\cap_{i=1}^{2m+1}\mathrm{Ker}\left(dx_{i}\right)\\
T^{-}Z^{+}&:= \cap_{i=2m+2}^{4m+2}
\mathrm{Ker}\left(dx_{i}\right).
\end{align*}
The para-complex structure ${\mathcal J}^{+}$ preserves the fibers
of the reflector projection $\pi^{+}:Z^{+}\rightarrow M$ and the
two distributions $\mathrm{Ker}(\pi^{+})_{*}\cap T^{+}Z^{+}$ and
$\mathrm{Ker}(\pi^{+})_{*}\cap T^{-}Z^{+}$ have rank one. Suppose
they are generated on $\mathcal U$ by two vector fields, say
$X^{+}$ and $X^{-}$ respectively. From our choice of $(x_{1},
\cdots , x_{n+2})$,
\begin{align*}
&dx_{1}(X^{+})=\cdots = dx_{2m+1}(X^{+})=0\\
&dx_{2m+2}(X^{-})=\cdots = dx_{4m+2}(X^{-})=0.\\
\end{align*}
Eventually changing the order of $(x_{1}, \cdots , x_{2m+1})$, we
suppose that $dx_{1}\left(X^{-}\right)\neq 0$ at $I_{p}$;
similarly, we can take $dx_{2m+2}\left(X^{+}\right)\neq 0$ at
$I_{p}.$ Let $\mathcal S$ be a codimension two submanifold of
$\mathcal U$, which contains $I_{p}$ and is defined by
\begin{align*}
x^{1}&= f(x^{2},\cdots , x^{2m+1})\\
x^{2m+2}&= g(x^{2m+3},\cdots , x^{4m+2}),
\end{align*}
where $f$, $g$ are smooth functions with all partial derivatives
$\frac{\partial f}{\partial x_{j}}$ ($2\leq j\leq 2m+1$) and
$\frac{\partial g}{\partial x_{j}}$ ($2m+3\leq j\leq 4m+2$) zero
at $I_{p}$. Then $\mathcal S$ intersects the fibers of $\pi^{+}$
transversally in a neighborhood of $I_{p}$ and the tangent bundle
$T{\mathcal S}$ is preserved by ${\mathcal J}^{+}$. It follows
that $\mathcal S$ is the image of a compatible almost para-complex
structure $I$ of $(M, {\mathcal P})$, viewed as a (local) section
of $\pi^{+}:Z^{+}\rightarrow M$. From Lemma \ref{preg1} the almost
para-complex structure $I$ is integrable. Clearly, $I$ extends
$I_{p}$ in a neighborhood of $p$. We proved that the first
condition implies the fourth. The fourth condition implies the
fifth. Our claim follows.
\end{proof}

\end{document}